\documentclass{gtmon_a}
\pdfoutput=1

\usepackage{mdwtab}

\proceedingstitle{The Zieschang Gedenkschrift}
\conferencestart{5 September 2007}
\conferenceend{8 September 2007}
\conferencename{Conference in honour of Heiner Zieschang}
\conferencelocation{Toulouse, France}

\editor{Michel Boileau}
\givenname{Michel}
\surname{Boileau}

\editor{Martin Scharlemann}
\givenname{Martin}
\surname{Scharlemann}

\editor{Richard Weidmann}
\givenname{Richard}
\surname{Weidmann}

\title{On multiplicity of mappings between surfaces}

\author{Semeon Bogatyi}
\givenname{Semeon}
\surname{Bogatyi}
\address{{\rm SB, EK:}\qua Department of Mathematics and Mechanics\\
Moscow State University\\\newline
Moscow 119992\\Russia\vspace{3pt}\\\newline
{\rm JF:}\qua Fachbereich 6 -- Mathematik\\
Universit\"at Siegen\\\newline
57068 Siegen\\Germany}
\email{bogatyi@inbox.ru}
\urladdr{}

\author{Jan Fricke}
\givenname{Jan}
\surname{Fricke}
\email{fricke@math.uni-siegen.de}
\urladdr{}

\author{Elena Kudryavtseva}
\givenname{Elena}
\surname{Kudryavtseva}
\email{eakudr@mech.math.msu.su}
\urladdr{}

\dedicatory{In grateful memory of Heiner, his wonderful collaboration and friendship}

\volumenumber{14}
\issuenumber{}
\publicationyear{2008}
\papernumber{03}
\startpage{49}
\endpage{62}

\doi{}
\MR{}
\Zbl{}

\arxivreference{}


\received{29 April 2006}
\revised{}
\accepted{9 November 2006}
\published{29 April 2008}
\publishedonline{29 April 2008}
\proposed{}
\seconded{}
\corresponding{}
\version{}

\makeatletter
\def\cnewtheorem#1[#2]#3{\newtheorem{#1}{#3}[section]
\expandafter\let\csname c@#1\endcsname\c@Thm}
\makeatother

\let\xysavmatrix\xymatrix
\def\xymatrix{\disablesubscriptcorrection\xysavmatrix}
\AtBeginDocument{\let\bar\wbar\let\tilde\wtilde\let\hat\what}
\AtBeginDocument{\renewcommand{\rho}{\varrho}\renewcommand{\Im}{\mathop{\mathrm{Im}}\nolimits}}
\newcommand{\IntD}{\smash{\mskip4mu\mathring{\mskip-4mu\vrule width0pt height7pt depth0pt\smash{D}}}}
\newcommand{\IntM}{\smash{\mskip1mu\mathring{\mskip-1mu\vrule width0pt height7pt depth0pt\smash{\M'}}}}
\AtBeginDocument{\def\N{N}}
\def\bGamma{\bar{\Gamma}}

\newtheorem{Thm}{Theorem}[section]
\cnewtheorem{Pro}[Thm]{Proposition}
\makeautorefname{Thm}{Theorem}
\makeautorefname{Pro}{Proposition}
\newcommand{\RR}{\mathbb R}
\newcommand{\CC}{\mathbb C}

\def\M{M} 
\def\n{n} 
\def\m{m} 
\def\MMR{\mathrm {MMR}}
\def\MI{\mathrm {MI}}
\def\NI{\mathrm {NI}}
\newcommand{\sgn}{\mathrm {sgn}}
\newcommand{\id}{\mathrm {id}}
\newcommand{\Proof} {\begin{proof}}

\renewcommand{\int} {{\mathrm {Int}}}

\begin{document}

\begin{asciiabstract}   
Let M and N be two closed (not necessarily orientable)
surfaces, and f a continuous map from M to N. By
definition, the minimal multiplicity MMR[f] of the map f
denotes the minimal integer k having the following property: f
can be deformed into a map g such that the number |g^{-1}(c)| of
preimages of any point c in N under g is at most k. We calculate
MMR[f] for any map $f$ of positive absolute degree A(f). The
answer is formulated in terms of A(f),
[pi_1(N):f_#(pi_1(M))], and the Euler characteristics of M
and N. For a map f with A(f)=0, we prove the inequalities
2 <= MMR[f] <= 4.
\end{asciiabstract}

\begin{htmlabstract}
Let M and N be two closed (not necessarily orientable) surfaces, and
f:M &rarr; N a continuous map. By definition, the <em>minimal
multiplicity</em> MMR[f] of the map f denotes the minimal integer k
having the following property: f can be deformed into a map g such
that the number &#x2223;g<sup>-1</sup>(c)&#x2223; of preimages of any
point c&isin; N under g is &le; k. We calculate MMR[f] for any map f
of positive absolute degree A(f). The answer is formulated in terms of
A(f), [&pi;<sub>1</sub>(N):f<sub>#</sub>(&pi;<sub>1</sub>(M))], and
the Euler characteristics of M and N. For a map f with A(f)=0, we
prove the inequalities 2&le;MMR[f]&le;4.
\end{htmlabstract}

\begin{abstract}   
Let $M$ and $N$ be two closed (not necessarily orientable)
surfaces, and $f \colon\thinspace   M \to N$ a continuous map. By
definition, the \textit{minimal multiplicity} $\mathrm{MMR}[f]$ of the map $f$
denotes the minimal integer $k$ having the following property: $f$
can be deformed into a map $g$ such that the number $|g^{-1}(c)|$ of
preimages of any point $c\in N$ under $g$ is $\le k$. We calculate
$\mathrm{MMR}[f]$ for any map $f$ of positive absolute degree $A(f)$. The
answer is formulated in terms of $A(f)$,
$[\pi_1(N):f_\#(\pi_1(M))]$, and the Euler characteristics of $M$
and $N$. For a map $f$ with $A(f)=0$, we prove the inequalities
$2\le\mathrm{MMR}[f]\le4$.
\end{abstract}

\maketitle

\section{Introduction} \label{sec:intro}

For a continuous map $f\co  X\to Y$ between topological spaces,
we define the \textit{multiplicity} of $f$ as \ $\max_{y\in Y}
|f^{-1}(y)|$, and the \textit{minimal multiplicity} of $f$ as the
minimal multiplicity of maps homotopic to $f$, that is
$$
\MMR[f]\,:=\,\min_{g\simeq f} \max_{y\in Y} |g^{-1}(y)|.
$$
From now on, $\simeq$ means that the mappings are homotopic. The
problem of determining $\MMR[f]$ arises. This problem is closely
related to the \textit{self-intersection problem} of determining the
\textit{minimal self-intersection number} (see Bogatyi, Kudryavtseva and Zieschang \cite{BKZ2,BKZ3})
 $$
\MI[f]\,:=\, \min_{g\simeq f} |\int(g)| , \quad
\int(g)\,:=\,\{(x,y)\in X\times X\,|\, x\ne y,\ g(x)=g(y)\} /
\Sigma_2
 $$
(here $\Sigma_2$ is the symmetric group in two symbols, which acts
on $X\times X$ by permutations of the coordinates), and to the
problem of determining the \textit{minimal (unordered) $\mu$--tuple
self-intersection number}
 $$
\MI_\mu[f]\,:=\, \min_{g\simeq f} |\int_\mu(g)| , \quad
\int_\mu(g)\,:=\,\{I\subset X \,|\, |I|=\mu,\ |g(I)|=1\}, \quad
\mu\ge 2.
 $$
Clearly, $\MI[f]=\MI_2[f]$, and one easily shows\footnote{(Indeed,
take a map $g\simeq f$ such that $\MI_\mu[f]=|\int_\mu(g)|=:\ell$.
We can assume that $\ell<\infty$. Then
$\smash{\ell=\sum_{\smash{i\ge\mu}}\sum_{\smash{y\in Y,\ |g^{-1}(y)|=i}} {i\choose\mu}}$.
Hence, for every nonvanishing summand in this sum, one has
$\smash{{i\choose\mu}\le\ell}$ and
$$\textstyle\smash{{i\choose\mu+1}={i\choose\mu}\frac{i-\mu}{\mu+1}
  <{i\choose\mu}\frac{i}{\mu} \le {i\choose\mu}^{2}\le\ell{i\choose\mu}}.$$
Therefore $\textstyle\MI_{\mu+1}[f]\le|\int_{\mu+1}(g)|
=\smash{\sum_{i>\mu}\sum_{y\in Y,\ |g^{-1}(y)|=i} {i\choose\mu+1}}$, 
which is at most $\ell
\smash{\sum_{i>\mu} \sum_{y\in Y,\ |g^{-1}(y)|=i} {i\choose\mu}\le
\ell^2}$.)} that $\MI_{\mu+1}[f]\le(\MI_\mu[f])^2$, $\mu\ge2$. 
The
connection between $\MMR[f]$ and $\MI_\mu[f]$ is illustrated by the
following properties:
 $$
\MI_\mu[f]=0 \ \iff \ \MMR[f]<\mu \qquad \mbox{and} \qquad
\MI_\mu[f]>0 \ \iff \ \MMR[f]\ge\mu.
 $$
In particular, $\MI[f]=0$ if and only if $\MMR[f]=1$. The numbers
$\MMR[f]$, $\MI[f]$, and $\MI_\mu[f]$, measure, in a sense,
``complexity'' of the self-intersection set $\int(f)$.

It is natural to consider the above problem for maps
$f\co  \M^\m\to \N^\n$ between closed connected (nonempty)
smooth manifolds, where $\m=\dim \M$, $\n=\dim \N$. The problem is
nontrivial for $0<\m\le \n\le 2\m$.

Hurewicz~\cite{Hur} proved that, if $X$ is an $\m$--dimensional
compact metric space and $\m+1\le \n \le 2\m$, then any continuous
map $f\co  X\to\RR^\n$ can be deformed, by means of an arbitrary
small perturbation, to a map $g\co  X\to\RR^\n$ of multiplicity
$\le [\frac{\n}{\n-\m}]$. A similar assertion is also valid if the
Euclidean space $\RR^\n$ is replaced by an arbitrary smooth manifold
$\N^\n$. Thus, for $\m<\n\le2\m$, we have
 \begin{equation}
 \refstepcounter{Thm} 
 \label {eq:codim>0}
\MMR[f]\le \left[ \frac{\n}{\n-\m} \right].
 \tag{\hbox{\bf\theThm}}
\end{equation}
This inequality follows by observing that, for a ``generic'' map
$g\co  \M\to \N$, the set $\int_{\mu+1}(g)\subset\M$ has
dimension $(\mu+1)\m-\mu \n$, which is negative (and, thus,
$\MMR[f]\le\mu$) if $\mu>\frac{\m}{\n-\m}$.

The special case $\n=2\m$ is the classical self-intersection problem
which gives rise to Whitney's work~\cite{Wh}. Here the
estimation~\eqref{eq:codim>0} gives $\MMR[f]\in\{1,2\}$, and
computing $\MMR[f]$ is equivalent to deciding whether $\MI[f]=0$,
ie\ whether the map $f$ is homotopic to an embedding. Namely, we
have $\MMR[f]=1$ if $\MI[f]=0$, and $\MMR[f]=2$ if $\MI[f]>0$. A
useful tool for deciding whether $\MI[f]=0$ is the \textit{Nielsen
self-intersection number} $\NI[f]$ of $f$~\cite{BKZ2,BKZ3}.
One can show by using the Whitney trick~\cite {Wh} that
$\MI[f]=\NI[f]$ if $\m\ge3$. But, if $m\le2$, one has only the
inequality $\MI[f]\ge\NI[f]$ (see~our papers with Zieschang \cite{BKZ2,BKZ3} for
$\m=1$). For $\m=1$, there are several combinatorial and geometric
methods for deciding whether a closed curve on a surface is
homotopic to a simple closed curve (see, for example,~Gon\c calves, Kudryavtseva and Zieschang \cite {GKZ3}
and references therein). An answer in terms of the Nielsen
self-intersection number is given in \fullref{thm:curve}. In
the remaining case $\m=2$, we only know
   that $\NI[f]>0$ implies $\MI[f]>0$ (and thus $\MMR[f]=2$), but the question
whether $\NI[f]=0$ implies $\MI[f]=0$ is still open.

The present paper studies the number $\MMR[f]$ mainly in the case
$\m=\n\le2$. Here $\MMR[f]$ is closely related to the \textit{absolute
degree} $A(f)$ (as defined in~Hopf~\cite{H} or~Epstein~\cite{Ep}; see
also Kneser~\cite{K}, Olum~\cite{Ol} and Skora~\cite{Sk}) of the map $f$. A definition of
the absolute degree is also given in Definition 3.7 in the paper by Gon\c calves, Kudryavtseva and Zieschang~\cite{gkz} of this
volume. \fullref{thm:circle}
computes the number $\MMR[f]$ for a self-mapping $f$ of a circle
($\m=\n=1$). In the case $\m=\n=2$ (mappings between closed
surfaces), the following results are obtained. We calculate
$\MMR[f]$ in terms of $A(f)$,
$\ell(f)\,:=\,[\pi_1(\N):f_\#\pi_1(\M)]$, and the Euler
characteristics of the surfaces, for any map $f\co  \M\to \N$
with $A(f)>0$ (\mbox{\fullref{thm:cover}} and \mbox{\fullref {thm:pinch}}). We
also estimate $\MMR[f]$ for any map $f$ with $A(f)=0$ 
(\mbox{\fullref{thm:A=0}}). In particular, we prove that
 \begin{align*}
\MMR[f]&\in\{A(f),A(f)+2\} &\mbox{if} \quad A(f)>0,
 \\
\MMR[f]&\in\{2,3,4\} &\mbox{if} \quad A(f)=0.
 \end{align*}
The authors do not know whether $\MMR[f]\ge A(f)$ if $\m=\n\ge3$.

\subsubsection*{Acknowledgements} This work was partially completed
during the visit of the first and the third authors in June--July
2005 at the Fakult\"at f\"ur Mathematik, Universit\"at Siegen,
Deutschland. The visits have been supported by the FIGS-project
``Niedrigdimensionale Topologie und geometrisch-topologische
Methoden in der Gruppentheorie'' (1st author), by the Grant of the
President RF, project NSh--4578.2006.1, and by a PIMS Fellowship
(3rd author).

\section[Computing MMR for mappings of a circle]{Computing $\MMR[f]$ for mappings of a circle} \label{sec:1}

Any map $f\co  S^1\to \N$ with $\dim \N\ge3$ is homotopic to an
embedding, thus $\MMR[f]=1$. Consider the cases $\dim \N=1,2$.

\begin{Thm} \label{thm:circle}
For any self-map $f\co  S^1\to S^1$,
$$
\MMR[f] = \left\{ \begin{array}{ll} |\deg f|, & \deg f\ne0, \\
                                    2,        & \deg f=0. \end{array} \right.
$$
\end{Thm}

\Proof We will identify the circle $S^1$ with the unit circle in the
complex plane $\CC$. Consider the projection $p\co  \RR\to S^1$,
$p(r)=e^{2\pi ir}$, $r\in\RR$, of the universal covering $\RR$ of
$S^1$ to $S^1$.

Suppose that \ $\deg f\ne0$. Then $f$ is homotopic to the map
sending $z\mapsto z^{\deg f}$, $z\in S^1$. Thus all points have
exactly $|\deg f|$ preimages, hence $\MMR[f]\le|\deg f|$. Let us
show that the number of preimages can not be reduced. Since \ $\deg
f\ne0$, for every point $s\in S^1$ there exists a point $t\in S^1$
such that $f(t)=s$. Let $r_0\in\RR$ be a point such that $p(r_0)=t$.
Consider a lifting $\tilde f\co  \RR\to\RR$ of $f\co  S^1\to
S^1$. Then $\tilde f(r_0+1)=\tilde f(r_0)+\deg f$, so by the
Intermediate Value Theorem, there exist points $r_1,\dots,r_{|\deg
f|-1}\in(r_0,r_0+1)$ such that $\tilde f(r_i)=\tilde
f(r_0)+j\,\sgn(\deg f)$, $1\le j\le |\deg f|-1$. Thus
$p(r_0),p(r_1),\dots,p(r_{|\deg f|-1})$ are different preimages of
$s$ under the mapping $f$. This shows $\MMR[f]\ge|\deg f|$.

Suppose that \ $\deg f=0$. Let us show that there exists $g\simeq f$
with $|g^{-1}(s)|\le2$ for any $s\in S^1$. Indeed, take $g$ to be
the map given by the following rule: $g(z)=z$ if $\Im z\ge0$,
$g(z)=\bar z$ if $\Im z\le0$. It remains to show that for any
$f\co  S^1\to S^1$, $\deg f=0$, there exists a point $s\in S^1$
with $|f^{-1}(s)|\ge2$. Such a map $f$ lifts to a map $\bar
f\co  S^1\to\RR$, thus it is enough to show that $\bar f$ is not
an embedding. This can be easily deduced by taking two points
$s_0,s_1\in S^1$ with $\bar f(s_0)=\min_{s\in S^1}\bar f(s)$,
$\bar f(s_1)=\max_{s\in S^1}\bar f(s)$, and applying the Intermediate
Value Theorem to the restriction of $\bar f$ to two segments in
$S^1$ having endpoints at $s_0,s_1$. \end{proof}

Consider a closed curve $f\co  S^1\to \N^2$ on a closed surface
$\N^2$. Then computing $\MMR[f]$ is equivalent to deciding whether
the homotopy class $[f]$ of the curve $f$ contains a simple closed
curve. Namely, $\MMR[f]=1$ if $[f]$ contains a simple curve, and
$\MMR[f]=2$ otherwise.

\begin{Thm}{\rm \cite{BKZ2,BKZ3}}\qua \label{thm:curve}
A closed curve $f\co  S^1\to \N^2$ on a closed surface $\N^2$ is
homotopic to a simple closed curve if and only if $\NI[f]=0$ and one
of the following conditions is fulfilled: the curve $f$ is not
homotopic to a proper power of any closed curve on $\N$, or $f\simeq
g^2$ for some orientation-reversing closed curve $g\co  S^1\to
\N$.\hfill
 \qedsymbol
\end{Thm}

An analogue of \fullref{thm:curve} was proved by Turaev and
Viro~\cite[Corollary~II]{TV}, in terms of the intersection index
introduced therein.

\section {$\MMR(f)$ for maps of positive degree between surfaces}
\label {sec:A>0}

In the following, $\M=\M^2$ and $\N=\N^2$ are arbitrary connected
closed surfaces, ie\ $2$--dimensional manifolds. By $\chi(\M)$, we
denote the Euler characteristic of $\M$. For a continuous mapping
$f\co  \M\to \N$, $A(f)$ denotes its \textit{absolute degree}
(see Hopf~\cite{H}, Epstein \cite{Ep}, Kneser \cite{K}, Olum \cite{Ol}, Skora \cite{Sk} or
Gon\c calves, Kudryavtseva  and Zieschang \cite{GKZ1}). Denote the index of the image of the fundamental group
of $\M$ in the fundamental group of $\N$ by
$\ell(f)\,:=\,[\pi_1(\N,f(x_0)):f_\#(\pi_1(\M,x_0))]$ for some
$x_0\in \M$. Actually the number $\ell(f)$ does not depend on the
choice of the point $x_0$.

The following consequence of Kneser's inequality will be central in
the proof of our main result.

\begin{Pro} \label{pro:Kneser}
If $f\co  \M\to \N$ has absolute degree $d=A(f)>0$ then there are
at most $d\cdot\chi(\N)-\chi(\M)$ points in $\N$ whose preimages
have cardinality $\le d-1$. Moreover, if pairwise different points
$y_1,\dots,y_r$ of $\N$ have $\mu_1,\dots,\mu_r$ preimages,
respectively, then
 $$
d\cdot\chi(\N)\ge\chi(\M)+\sum_{i=1}^r (d-\mu_i).
 $$
\end{Pro}

\Proof In the case when $r=1$ and $f$ is orientation-true, the
latter inequality was proved in Theorem~2.5~(a) of~\cite{GKZ1}. In
the general case, the inequality can be proved using the techniques
in~\cite{BGKZ,GKZ1,GZ,BGZ}, as follows.

If $f$ is not orientation-true and $d=A(f)>0$ then $d=\ell(f)$, due
to the result of Kneser~\cite{K1928,K}. On the other hand,
one has $\mu_i\ge\ell(f)$, $1\le i\le r$, since the map $f$ admits a
lifting $\hat f\co  M\to\hat N$ such that $f=p\circ\hat f$, where
$p\co  \hat N\to N$ is an $\ell(f)$--fold covering corresponding
to the subgroup $f_\#(\pi_1(\M,x_0))$ of $\pi_1(\N,f(x_0))$, and
$A(\hat f)=1$, hence $\hat f$ is surjective. Therefore
$\sum_{i=1}^r(d-\mu_i)\le0$. This, together with the Kneser
inequality~\cite {K}, $d\cdot\chi(\N)\ge\chi(\M)$, implies the
desired inequality.

If $f$ is orientation-true, one proceeds as in the proof of
Proposition~2.5~(a) of~\cite {GKZ1}, where one replaces the single
point $y_0\in\N$ by the set of $r$ points $y_1,\dots,y_r$. More
specifically, by applying a suitable deformation, one can assume
that there are small pairwise disjoint disks $D_i,D_{ij}$, $1\le
i\le r$, $1\le j\le \mu_i$, around the points $y_i$ of $\N$ and the
points of $f^{-1}(y_i)$ such that $\smash{f^{-1}(\IntD_i) =
\bigcup_{j=1}^{\mu_i} \IntD_{ij}}$, and $f|_{D_{ij}}$ is a branched
covering of type $\smash{z\mapsto z^{d_{ij}}}$ for some positive integer
$d_{ij}$. Therefore the complement of these open disks are two
compact surfaces $F\subset \M$, $G\subset \N$ such that the
restriction of $f$ induces a proper map carrying the boundary into
the boundary, $f|_{F}\co   (F, \partial F) \to (G, \partial G)$.
By Proposition~1.6 of~\cite{GKZ1} (or by a more general Theorem~4.1
of~\cite{Sk}), $\chi(F) \leq A(f)\cdot \chi(G)$. This, together with
$\chi(F)=\chi(\M)-\sum_{i=1}^r\mu_i$, $\chi(G)=\chi(\N)-r$, gives
the desired inequality. \end{proof}

\begin{Thm} \label{thm:cover}
Suppose that $f\co  \M\to \N$ has absolute degree \ $d=A(f)>0$.
If \ $\ell(f)\ne d$, or \ $\ell(f)=d$ \ and \
$d\cdot\chi(\N)=\chi(\M)$, then \ $\MMR[f]=d$.
\end{Thm}

\Proof The inequality $\MMR[f]\ge A(f)$ follows from the first part
of \fullref {pro:Kneser}.

Let us show the converse inequality, $\MMR[f]\le A(f)$. It follows
from~\cite {Ed,Sk,K}, respectively, that the
mapping $f$ is homotopic to a $d$--fold covering which is branched in
the first case and unbranched in the second case. Thus, we found a
mapping which is homotopic to $f$, and the preimage of any point of
$N$ has cardinality $\le d$. \end{proof}

\begin{Thm} \label{thm:pinch}
Suppose that $f\co  \M\to \N$ has absolute degree \ $d=A(f)>0$.
If \ $\ell(f)=d$ \ and \ $d\cdot\chi(\N)\ne\chi(\M)$, then \
$\MMR[f]=d+2$.
\end{Thm}

\Proof \textbf{Case 1}\qua Suppose that $d=A(f)=1$. It follows from~\cite {Ed,Sk}
that the mapping $f$ is homotopic to a pinching map where
the pinched subsurface $\M'\subset \M$, $\partial \M'\simeq S^1$, is
different from the $2$--disk $D^2$ (here the natural projection $\M\to
\M/\M'$ is called a pinching map).

\begin{figure}[ht!]
\setlength{\unitlength}{9pt}
\begin{center}
\begin{picture} (10,11.5)(-5,-1.5) 
\small
\put(-8,-2){
\put(0,-5.4){ \thicklines \qbezier[80](-3.7,8.9)(-3.3,9)(-2.85,9.09)
\qbezier[300](-1.55,9.25)(0,9.45)(1.55,9.25)
\qbezier[80](3.7,8.9)(3.3,9)(2.85,9.09)

\qbezier[400](-3.7,8.9)(-6.3,8)(-3.7,7.1)
\qbezier[400](3.7,8.9)(6.3,8)(3.7,7.1)
\qbezier[400](-3.7,7.1)(0,6.2)(3.7,7.1) }
\thinlines \put(-2.2,1){ \qbezier[200](-1,4)(-.3,2.7)(-.88,1.8)
\qbezier[200](-1,4)(-2.32,7)(-1.28,10)
\qbezier[200](-1.28,10)(0,12)(1.28,10)
\qbezier[200](1.28,10)(2.32,7)(1,4)
\qbezier[200](1,4)(.3,2.7)(.88,1.8) }
\put(-2.2,3.5){ \qbezier[50](-.8,7)(0,6.2)(.8,7)
\qbezier[50](-.5,6.8)(0,7.2)(.5,6.8) } \put(-2.2,1.6){
\qbezier[50](-.8,7)(0,6.2)(.8,7)
\qbezier[50](-.5,6.8)(0,7.2)(.5,6.8) } \put(-2.2,-.3){
\qbezier[50](-.8,7)(0,6.2)(.8,7)
\qbezier[50](-.5,6.8)(0,7.2)(.5,6.8) }
\put(2.2,1){ \qbezier[200](-1,4)(-.3,2.7)(-.88,1.8)
\qbezier[200](-1,4)(-2.32,7)(0,9) \qbezier[200](0,9)(2.32,7)(1,4)
\qbezier[200](1,4)(.3,2.7)(.88,1.8)
\put(0,9) {\line(0,-1){2.9}} \put(0,6.1){
\qbezier[30](0,0)(-.4,-.35)(-.8,-.4) 
\qbezier[30](-.8,-.4)(-1.4,-.4)(-1.48,0)
\qbezier[7](0,0)(-.5,.35)(-.8,.4)
\qbezier[8](-.8,.4)(-1.4,.4)(-1.48,0)
\qbezier[30](0,0)(.4,-.35)(.8,-.4) 
\qbezier[30](.8,-.4)(1.4,-.4)(1.48,0)
\qbezier[7](0,0)(.5,.35)(.8,.4) \qbezier[8](.8,.4)(1.4,.4)(1.48,0) }
\put(0,7.5){
\qbezier[30](0,0)(-.3,-.25)(-.6,-.3) 
\qbezier[30](-.6,-.3)(-1.06,-.3)(-1.14,0)
\qbezier[5](0,0)(-.3,.25)(-.6,.3)
\qbezier[6](-.6,.3)(-1.06,.3)(-1.14,0)
\qbezier[30](0,0)(.3,-.25)(.6,-.3) 
\qbezier[30](.6,-.3)(1.06,-.3)(1.14,0)
\qbezier[5](0,0)(.3,.25)(.6,.3) \qbezier[6](.6,.3)(1.06,.3)(1.14,0)
} } } \put(-1,4){$\longrightarrow$}
\put(8,4) 
{ \thicklines \qbezier[400](-3.1,3.1)(0,5.7)(3.1,3.1)
\qbezier[300](-3.1,3.1)(-5.7,0)(-3.1,-3.1)
\qbezier[300](3.1,3.1)(5.7,0)(3.1,-3.1)
\qbezier[400](-3.1,-3.1)(0,-5.7)(3.1,-3.1) \put(.3,-5.2){
\thinlines \put(-1.7,1){ \qbezier[200](-.6,2.7)(-.2,1.8)(-.4,1.2)
\qbezier[200](-.6,2.7)(-1.54,4.7)(-.84,6.66)
\qbezier[200](-.84,6.66)(0,8)(.84,6.66)
\qbezier[200](.84,6.66)(1.54,4.7)(.6,2.7)
\qbezier[200](.6,2.7)(.2,1.8)(.4,1.2)
\qbezier[5](-.4,1.2)(0,1.45)(.4,1.2) }
\put(-1.7,7.3){\circle{.6}} \put(-1.7,6){\circle{.6}}
\put(-1.7,4.7){\circle{.6}}
\put(1.7,1.5){ \qbezier[200](-.66,2.7)(-.2,1.8)(-.4,1.2)
\qbezier[200](-.66,2.7)(-1.54,4.7)(0,6)
\qbezier[200](0,6)(1.54,4.7)(.66,2.7)
\qbezier[200](.66,2.7)(.2,1.8)(.4,1.2)
\qbezier[5](-.4,1.2)(0,1.4)(.4,1.2)
\put(0,6) {\line(0,-1){2}} } } }
\end{picture}
\end{center}
\label{figure1}
\caption{A proper map \ $g'\co \M'\to D^2$ \ of
multiplicity 3}
\end{figure}
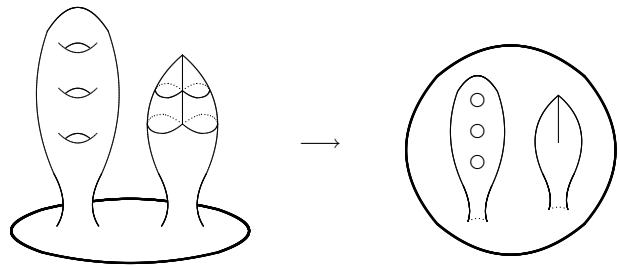

Let us show that such a pinching map is homotopic to a map $g$ of
multiplicity $\le3$. For this, we construct a proper continuous map
$g'\co  (\M',\partial \M')\to (D^2,\partial D^2)$ whose
restriction to $\partial \M'$ is a homeomorphism, and whose
multiplicity equals $3$. Such a map $g'$ is shown in \fullref{figure1}. We may
identify $\N$ with the surface which is obtained by gluing of
$\M\setminus \IntM$ and $D^2$ by means of the aforementioned
homeomorphism of the boundary circles, where $\IntM$ denotes the
interior of $\M'$. Define $g\co  \M\to \N$ as $g|_{\M\setminus
\M'}=\id_{\M\setminus \M'}$ and $g|_{\M'}=g'$. Clearly, $f\simeq g$,
since $g'$ is homotopic relative boundary to a pinching map. In
Case~2 below, we will use the following property of the constructed
map $g$: its restriction to the preimage of the complement
$\N\setminus D^2$ of a disk is injective.

It follows from the inequality of Euler characteristics of $\M$ and
$\N$ that $f$ is not homotopic to an embedding. (Indeed, otherwise
such an embedding $g$ is a homeomorphism onto $g(M)$; it follows from Brouwer's
Theorem on Invariance of Domain~\cite{Brouwer} that $g$ is
surjective and, therefore, it is a homeomorphism.) Suppose that $f$
is homotopic to a map $g\co  \M\to \N$ of multiplicity 2, we will
show that this leads to a contradiction. Let $y\in \N$ be a point
with $g^{-1}(y)=\{x_1,x_2\}$. Then the local degree of $g$ at each
of the points $x_1$ and $x_2$ is defined modulo~2, and
 $$
\deg(g,x_1)+\deg(g,x_2)\equiv A(g)\equiv A(f)\equiv1 \mod 2.
 $$
Without loss of generality, we may assume that $\deg(g,x_1)\ne0$.
This implies that the image of any neighbourhood of $x_1$ contains a
neighbourhood of $y=g(x_1)$, since otherwise one could construct a
map $F\co  D^2\to S^1$ with $\deg(F|_{\partial
D^2})=\deg(g,x_1)\ne0$. Therefore the restriction of $g$ to an
appropriate neighbourhood of $x_2$ is injective and, thus (by
Brouwer's Theorem on Invariance of Domain~\cite{Brouwer}), is a
homeomorphism onto a neighbourhood of $y$. This implies that
$\deg(g,x_2)=\pm1$. Similar arguments show that $\deg(g,x_1)=\pm1$,
a contradiction.

\textbf{Case 2}\qua Suppose that $d=A(f)=\ell(f)\ge2$. Let us construct a map
$g$ which is homotopic to $f$ and has multiplicity $A(f)+2$.
Consider a covering $p\co  \tilde \N\to \N$ which corresponds to
the subgroup $f_\#(\pi_1(\M,x_0))$ of $\pi_1(\N,f(x_0))$. So, this
is an $\ell(f)$--fold covering. Let $y\in \N$ be an arbitrary point
and $D$ a small closed neighbourhood which is homeomorphic to the
disk $D^2$. Let $D_1,\dots,D_d$ be the connected components of
$p^{-1}(D)$.

Let $\tilde f\co  \M\to\tilde \N$ be a lifting of $f$. Then
$A(\tilde f)=\ell(\tilde f)=1$. By Case~1, there exists a map
$\tilde g\co  \M\to\tilde \N$ which is homotopic to $\tilde f$
and has multiplicity $\le3$. Then the map $g\,:=\,p\circ\tilde g$ is
homotopic to $f=p\circ\tilde f$. By Case~1, we may also assume that
$\tilde g$ is injective on $\tilde g^{-1}(\tilde \N \setminus D_1)$.
Therefore the map $g$ has multiplicity $\ell(f)+2=A(f)+2$.

Let us show that the multiplicity of $f$ is $\ge\ell(f)+2$. Let
$\tilde f\co  \M\to\tilde \N$ be a lifting of $f$ to this
$\ell(f)$--fold covering, thus $A(\tilde f)=\ell(\tilde f)=1$. By
Case~1, there exists a point $\tilde y\in\tilde \N$ whose preimage
under $\tilde f$ has cardinality $\ge3$. Since $A(\tilde f)>0$,
every point of $p^{-1}(p(\tilde y))$ has a nonempty preimage under
$\tilde f$. Therefore $f^{-1}(p(\tilde y))$ has cardinality
at least $\ell(f)+2=A(f)+2$. \end{proof}

\section[Estimates for MMR(f) if A(f)=0]{Estimates for $\MMR(f)$ if $A(f)=0$} \label{sec:A=0}

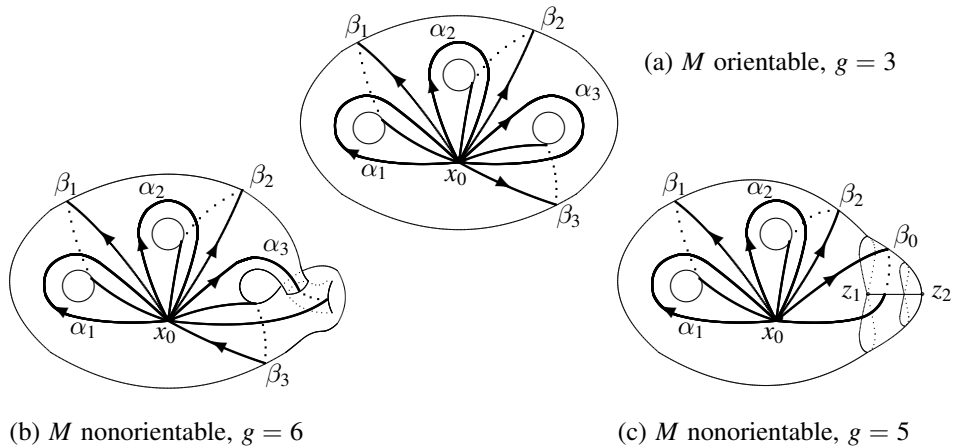
\begin{figure}[ht!]
\setlength{\unitlength}{10pt}
\begin{center}
\begin{picture} (10,17)(-5,-12.5) 
\small
\put(-1,0){
\qbezier[300](-4.44,2.7)(0,5.4)(4.44,2.7) 
\qbezier[200](-4.44,2.7)(-7.56,0)(-4.44,-2.7)
\qbezier[200](4.44,2.7)(7.56,0)(4.44,-2.7)
\qbezier[300](-4.44,-2.7)(0,-5.4)(4.44,-2.7)

\put(-3.4,-.2){\circle{1.2}} 
\put(0,1.8){\circle{1.2}} \put(3.4,-.2){\circle{1.2}}

\put(0,-1.5){\put(-.18,-.18){\small$\bullet$}} 
\put(-.6,-2.2){\small$x_0$} \thicklines
\put(0,-1.5){ 
\qbezier[200](0,0)(-2.2,-.2)(-3.7,.2)
\put(-3.8,-.5){\small$\alpha_1$}
\qbezier[200](-3.7,.2)(-4.4,.4)(-4.65,1)
\put(-4.2,.45){\vector(-2,1){.2}}
\qbezier[200](-4.65,1)(-4.9,2.2)(-3.7,2.5)
\qbezier[200](-3.7,2.5)(-3.2,2.6)(-2.6,2.2)
\qbezier[200](0,0)(-1.4,1.3)(-2.6,2.2) }
\put(0,-1.5){ 
\qbezier[100](0,0)(-2.2,3.3)(-3.8,4.5)
\put(-4.4,4.9){\small$\beta_1$}
\qbezier[8](-3.8,4.5)(-3.6,3)(-2.9,1.6)
\put(-2.37,3.2){\vector(-1,1){.2}}
\qbezier[100](0,0)(-1.6,.5)(-2.9,1.6) }
\put(0,-1.5){ 
\qbezier[200](0,0)(-.8,1.5)(-1.1,3.1)
\put(-.9,2.45){\vector(-1,4){.2}}
\qbezier[200](-1.1,3.1)(-1.2,3.7)(-.8,4.2)
\put(-1.2,4.8){\small$\alpha_2$}
\qbezier[200](-.8,4.2)(0,4.9)(.8,4.2)
\qbezier[200](1.1,3.1)(1.2,3.7)(.8,4.2)
\qbezier[200](0,0)(.8,1.5)(1.1,3.1) }
\put(0,-1.5){ 
\qbezier[100](0,0)(1.8,2.5)(2.8,4.95) \put(3,5.2){\small$\beta_2$}
\qbezier[8](2.8,4.95)(1.3,4.3)(.5,3) \put(1.7,2.6){\vector(1,2){.2}}
\qbezier[100](0,0)(.2,1.6)(.5,3) }
\put(0,-1.5){ 
\qbezier[200](0,0)(2.2,-.2)(3.7,.2)
\qbezier[200](3.7,.2)(4.4,.4)(4.65,1) \put(4.3,2.5){\small$\alpha_3$}
\qbezier[200](4.65,1)(4.9,2.2)(3.7,2.5)
\qbezier[200](3.7,2.5)(3.2,2.6)(2.6,2.2)
\qbezier[200](0,0)(1.4,1.3)(2.6,2.2)
\put(1.85,1.57){\vector(1,1){.2}} }
\put(0,-1.5){ 
\qbezier[100](0,0)(1.7,.7)(3.3,.65)
\qbezier[8](3.3,.65)(3.8,-.4)(3.65,-1.6)
\put(1.8,-.87){\vector(2,-1){.2}}
\qbezier[100](0,0)(1.8,-1)(3.65,-1.6) \put(3.6,-2.3){\small$\beta_3$}
} 
   \put(7,2){\small (a) $\M$ orientable, $g=3$} }

\put(-12,-6){
\qbezier[300](-4.44,2.7)(-.6,4.95)(2.8,3.45) 
\qbezier[200](-4.44,2.7)(-7.56,0)(-4.44,-2.7) 
\qbezier[300](-4.44,-2.7)(0,-5.4)(4.44,-2.7) 

\put(-3.4,-.2){\circle{1.2}} 
\put(0,1.8){\circle{1.2}}
\qbezier[100](3.3,-.85)(3.8,-.85)(4.25,-.2) 
\qbezier[6](4.25,-.2)(4.7,.2)(5.1,.4)  

\qbezier[100](3.3,-.85)(2.75,-.75)(2.7,-.2) 
\qbezier[100](3.3,.45)(2.75,.35)(2.7,-.2) 
\qbezier[100](3.3,.45)(4,.5)(4.5,-.5) 
\qbezier[10](4.5,-.5)(5,-1.2)(6,-1.2) 

\qbezier[100](6.2,0)(5.9,-.6)(6.2,-1.2) 
\qbezier[50](5.1,.4)(5.9,.76)(6.4,.4) 
\qbezier[50](6.4,.4)(7,-.6)(6.4,-1.6) 
\qbezier[100](5,.6)(4.7,2.3)(2.8,3.45) 
\qbezier[50](5,.6)(5.1,.2)(5.3,0) 
\qbezier[5](5.3,0)(5.65,-.3)(6,0) 

\qbezier[50](5,-2.25)(4.7,-2.5)(4.44,-2.7) 
\qbezier[50](5,-2.25)(5.4,-1.9)(5.8,-1.9) 
\qbezier[50](6.4,-1.6)(6.1,-1.9)(5.8,-1.9)

\qbezier[50](4.5,-.5)(5.1,-.5)(5.3,0) 
\qbezier[6](4.5,-.5)(4.7,0)(5.3,0) 

\put(0,-1.5){\put(-.18,-.18){\small$\bullet$}} 
\put(-.6,-2.2){\small$x_0$} \thicklines
\put(0,-1.5){ 
\qbezier[200](0,0)(-2.2,-.2)(-3.7,.2)
\put(-3.8,-.5){\small$\alpha_1$}
\qbezier[200](-3.7,.2)(-4.4,.4)(-4.65,1)
\put(-4.2,.45){\vector(-2,1){.2}}
\qbezier[200](-4.65,1)(-4.9,2.2)(-3.7,2.5)
\qbezier[200](-3.7,2.5)(-3.2,2.6)(-2.6,2.2)
\qbezier[200](0,0)(-1.4,1.3)(-2.6,2.2) }
\put(0,-1.5){ 
\qbezier[100](0,0)(-2.2,3.3)(-3.8,4.5)
\put(-4.4,4.9){\small$\beta_1$}
\qbezier[8](-3.8,4.5)(-3.6,3)(-2.9,1.6)
\put(-2.37,3.2){\vector(-1,1){.2}}
\qbezier[100](0,0)(-1.6,.5)(-2.9,1.6) }
\put(0,-1.5){ 
\qbezier[200](0,0)(-.8,1.5)(-1.1,3.1)
\put(-.9,2.45){\vector(-1,4){.2}}
\qbezier[200](-1.1,3.1)(-1.2,3.7)(-.8,4.2)
\put(-1.2,4.8){\small$\alpha_2$}
\qbezier[200](-.8,4.2)(0,4.9)(.8,4.2)
\qbezier[200](1.1,3.1)(1.2,3.7)(.8,4.2)
\qbezier[200](0,0)(.8,1.5)(1.1,3.1) }
\put(0,-1.5){ 
\qbezier[100](0,0)(1.8,2.5)(2.8,4.95) \put(3,5.2){\small$\beta_2$}
\qbezier[8](2.8,4.95)(1.3,4.3)(.5,3) \put(1.7,2.6){\vector(1,2){.2}}
\qbezier[100](0,0)(.2,1.6)(.5,3) }
\put(0,-1.5){ 
\qbezier[200](0,0)(4.2,-.5)(6,.8) \put(3.5,2.65){\small$\alpha_3$}
\qbezier[100](4.95,1.15)(4.4,2.3)(3.3,2.4)
\qbezier[4](6,.8)(5.3,.8)(4.95,1.15)

\qbezier[200](3.3,2.4)(3.1,2.4)(2.6,2.2)
\qbezier[200](0,0)(1.4,1.3)(2.6,2.2)
\put(1.85,1.57){\vector(1,1){.2}} }
\put(0,-1.5){ 
\qbezier[100](0,0)(1.7,.7)(3.3,.65)
\qbezier[8](3.3,.65)(3.8,-.4)(3.65,-1.6)
\put(1.8,-.87){\vector(-2,1){.2}}
\qbezier[100](0,0)(1.8,-1)(3.65,-1.6) \put(3.6,-2.3){\small$\beta_3$}
} \put(-6,-6){\small (b) $\M$ nonorientable, $g=6$} }

\put(11,-6){
\qbezier[300](-4.44,2.7)(0,5.4)(3.44,1.7) 
\qbezier[200](-4.44,2.7)(-7.56,0)(-4.44,-2.7) 
\qbezier[200](3.44,1.7)(7.56,-.5)(3.44,-2.7) 
\qbezier[300](-4.44,-2.7)(-.3,-5.2)(3.44,-2.7) 

\put(-3.4,-.2){\circle{1.2}} 
\put(0,1.8){\circle{1.2}}
\put(3.44,-.5){\line(1,0){2.1}} \put(5.53,-.5){\circle{.06}}
\put(5.8,-.6){\small$z_2$} \put(3.47,-.5){\circle{.06}}
\put(2.4,-.6){\small$z_1$}
\put(3.44,-.5){ 
\qbezier[30](0,0)(-.2,.75)(-.3,1.5) 
\qbezier[30](-.3,1.5)(-.3,2.1)(0,2.2)
\qbezier[7](0,0)(.2,.75)(.3,1.5) \qbezier[8](.3,1.5)(.3,2.1)(0,2.2)
\qbezier[30](0,0)(-.2,-.75)(-.3,-1.5) 
\qbezier[30](-.3,-1.5)(-.3,-2.1)(0,-2.2)
\qbezier[7](0,0)(.2,-.75)(.3,-1.5)
\qbezier[5](.3,-1.5)(.3,-2.1)(0,-2.2) }
\put(4.84,-.5){ 
\qbezier[20](0,0)(-.1,.425)(-.17,.85) 
\qbezier[10](-.17,.85)(-.17,1.15)(0,1.2)
\qbezier[5](0,0)(.1,.425)(.17,.85)
\qbezier[3](.17,.85)(.17,1.15)(0,1.2)
\qbezier[20](0,0)(-.1,-.425)(-.17,-.85) 
\qbezier[10](-.17,-.85)(-.17,-1.15)(0,-1.2)
\qbezier[5](0,0)(.1,-.425)(.17,-.85)
\qbezier[3](.17,-.85)(.17,-1.15)(0,-1.2) }
\put(0,-1.5){\put(-.18,-.18){\small$\bullet$}} 
\put(-.6,-2.2){\small$x_0$}
\thicklines 
\put(0,-1.5){ 
\qbezier[200](0,0)(-2.2,-.2)(-3.7,.2)
\put(-3.8,-.5){\small$\alpha_1$}
\qbezier[200](-3.7,.2)(-4.4,.4)(-4.65,1)
\put(-4.2,.45){\vector(-2,1){.2}}
\qbezier[200](-4.65,1)(-4.9,2.2)(-3.7,2.5)
\qbezier[200](-3.7,2.5)(-3.2,2.6)(-2.6,2.2)
\qbezier[200](0,0)(-1.4,1.3)(-2.6,2.2) }
\put(0,-1.5){ 
\qbezier[100](0,0)(-2.2,3.3)(-3.8,4.5)
\put(-4.4,4.9){\small$\beta_1$}
\qbezier[8](-3.8,4.5)(-3.6,3)(-2.9,1.6)
\put(-2.37,3.2){\vector(-1,1){.2}}
\qbezier[100](0,0)(-1.6,.5)(-2.9,1.6) }
\put(0,-1.5){ 
\qbezier[200](0,0)(-.8,1.5)(-1.1,3.1)
\put(-.9,2.45){\vector(-1,4){.2}}
\qbezier[200](-1.1,3.1)(-1.2,3.7)(-.8,4.2)
\put(-1.2,4.8){\small$\alpha_2$}
\qbezier[200](-.8,4.2)(0,4.9)(.8,4.2)
\qbezier[200](1.1,3.1)(1.2,3.7)(.8,4.2)
\qbezier[200](0,0)(.8,1.5)(1.1,3.1) }
\put(0,-1.5){ 
\qbezier[100](0,0)(1.8,2.5)(2.35,4.15) \put(2.3,4.5){\small$\beta_2$}
\qbezier[6](2.35,4.15)(1.2,4.3)(.5,3)
\put(1.7,2.6){\vector(1,2){.2}} \qbezier[100](0,0)(.2,1.6)(.5,3) }
\put(0,-1.5){ 
\qbezier[200](0,0)(2.2,-.2)(3.3,.15)
\qbezier[80](3.3,.15)(3.9,.4)(4.1,1) \put(4.35,2.95){\small$\beta_0$}
\qbezier[5](4.1,1)(4.4,2)(4.2,2.7)
\qbezier[200](0,0)(3,2.5)(4.2,2.7) \put(2.2,1.65){\vector(3,2){.2}}
} \put(-6,-6){\small (c) $\M$ nonorientable, $g=5$} }
\end{picture}
\end{center}
\caption{A canonical system of cuts on a closed surface
$\M$}\label{figure2}
\end{figure}

Suppose that $\M$ is a connected orientable closed surface of genus
$g\ge0$. Consider the standard presentation of the closed surface
$\M$ as the boundary of a solid surface in $\RR^3$ which is obtained
from a closed 3-ball by attaching $g$ solid handles; see \mbox{\fullref{figure2}~(a)}.
Choose a base point $x_0\in \M$ and consider a system of simple
closed curves $\alpha_1,\beta_1,\dots,\alpha_g,\beta_g$ on $\M$
based at $x_0$ which form a \textit{canonical system of cuts}; see
\fullref{figure2}~(a). Then the fundamental group $\pi_1(\M,x_0)$ admits a
canonical presentation
 $$
\pi_1(\M,x_0)=\bigg\langle a_1,b_1,\dots,a_g,b_g\ \bigg|\
\prod_{j=1}^g[a_j,b_j] \bigg\rangle ,
 $$
where $a_j,b_j$ are the homotopy classes of the based loops
$\alpha_j,\beta_j$, respectively. Denote by $V_g$ the bouquet of $g$
circles $\alpha_1\cup\ldots\cup\alpha_g$ if $g\ge1$, $V_0:=\{x_0\}$
if $g=0$, and by $\rho$ a retraction $\rho\co  \M\to V_g$ which
maps all loops $\beta_j$ to the point $x_0$. We can assume that the
curves $\alpha_1,\dots,\alpha_g$ are contained in the plane
$\Pi\subset\RR^3$ which is tangent to $\M$ at $x_0$. (In \fullref{figure2}, the
plane $\Pi$ is parallel to the plane of the picture.)

Let $i\co  M\to\RR^3$ denote the inclusion, and
$p_\Pi\co  \RR^3\to\Pi$ the orthogonal projection. The following
properties of the map $p=p_\Pi\circ i\co  \M\to\Pi$ can be
assumed without loss of generality, and will be used later:

(p1)\qua The restriction of $p$ to a neighbourhood $U$ of the
base point $x_0\in \M$ is a homeomorphism onto a neighbourhood of
the point $p(x_0)$ in $\Pi$. Moreover, $p|_{V_r}\co  V_r\to\Pi$
is an embedding, and all curves $p|_{\alpha_j}\co  \alpha_j\to
\Pi$ are regular;

(p2)\qua All curves $p|_{\beta_j}$ are contractible in
$p(\M)$;

(p3)\qua $p(\M)$ is a regular neighbourhood of the graph
$p(V_r)$ in $\Pi$;

(p4)\qua The map $p$ has multiplicity $2$.

\eject
Suppose that $\M$ is a connected nonorientable closed surface of
genus $g\ge1$. Choose a base point $x_0\in \M$. Then the fundamental
group of $\M$ admits the following canonical presentation:
\begin{align*}
\pi_1(\M,x_0)&=\bigg\langle a_1,b_1,\dots,a_{\unfrac g2},b_{\unfrac g2}\
\bigg|\ \bigg(\prod_{j=1}^{\unfrac g2-1}[a_j,b_j]\bigg) \cdot
[a_{\unfrac g2},b_{\unfrac g2}]_- \bigg\rangle
 &&\mbox{if $g$ is even,}
\\
\pi_1(\M,x_0)&=\bigg\langle a_1,b_1,\dots,a_{[\unfrac g2]},b_{[\unfrac g2]},b_0\ \bigg|\ \bigg(\prod_{j=1}^{\upnfrac {g-1}2}[a_j,b_j]\bigg)
\cdot b_0^2 \bigg\rangle
 &&\mbox{if $g$ is odd,}
\end{align*}
where we use the notation
 $$
[x,y]\,=\,xyx^{-1}y^{-1}, \qquad [x,y]_-\,=\,xyx^{-1}y .
 $$
This presentation of the group $\pi_1(\M,x_0)$ corresponds to a
system of simple closed curves
$\alpha_1,\beta_1,\dots,\alpha_{[g/2]},\beta_{[g/2]},\beta_0$ on
$\M$ based at $x_0$, which form a \textit{canonical system of cuts};
see \fullref{figure2}~(b),~(c). Here the curve $\beta_0$ appears only if $g$ is
odd. Denote by $V_r$ the bouquet of $r=[g/2]$ circles
$\alpha_1\cup\ldots\cup\alpha_{[g/2]}$ for $g\ge2$, $V_0=\{x_0\}$
for $g=1$, and by $\rho$ a retraction $\rho\co  \M\to V_r$ which
maps all loops $\beta_j$ to the point $x_0$. We consider a
realization of $\M$ in $\RR^3$ via a map $i\co  \M\to\RR^3$ which
is an immersion if $g$ is even (see \fullref{figure2}~(b)), while, for $g$ odd,
the restriction $i|_{\M\setminus\{z_1,z_2\}}$ to the complement of
the set of two points $z_1,z_2\in\M\setminus\{x_0\}$ is an
immersion; see \fullref{figure2}~(c). We can assume that $i|_{V_r}$ is an
embedding with $i(V_r)\subset\Pi$, moreover $\Pi$ coincides with the
tangent plane to $i(\M)$ at $i(x_0)$.

\begin{figure}[ht!]
\setlength{\unitlength}{10pt}
\begin{center}
\begin{picture} (10,11)(-5,-6.5) 
\small
\put(-12.5,0){
\qbezier[300](-4.44,2.7)(0,5.4)(4.44,2.7) 
\qbezier[200](-4.44,2.7)(-7.56,0)(-4.44,-2.7) 
\qbezier[300](-4.44,-2.7)(0,-5.4)(4.44,-2.7) 
\put(3.8,-3.8){$p(\M)$}

\put(-3.4,-.2){\circle{1.2}} 
\put(0,1.8){\circle{1.2}}
\qbezier[100](3.3,-.85)(3.8,-.85)(4.25,-.2) 
\qbezier[100](3.3,-.85)(2.75,-.75)(2.7,-.2) 
\qbezier[100](3.3,.45)(2.75,.35)(2.7,-.2) 
\qbezier[100](3.3,.45)(3.8,.45)(4.25,-.2) 
\qbezier[50](5.7,1.2)(6,.025)(5.7,-1.15) 
\qbezier[50](5.7,1.2)(5.3,2.13)(4.44,2.7) 
\qbezier[50](5.7,-1.15)(5.3,-2.105)(4.44,-2.7) 

\put(0,-1.5){\put(-.18,-.18){\tiny$\bullet$}} 
\put(-.9,-2.4){$p(x_0)$}

\put(4.32,-.2) { \put(-.4,.4){\line(2,1){1.85}}
\put(-.4,-.4){\line(4,-1){1.85}} \put(2.05,1.325){\line(0,-1){2.2}}
\put(2.3,.1){$T$} \put(2.05,1.3){\line(-1,0){.2}}
\put(2.05,-.86){\line(-1,0){.2}} } \thicklines \put(4.25,-.2)
{\line(5,1){1.6}} \put(4.8,-.05){\circle{.06}}
\put(3.15,-2.2){$p(t)$}\thinlines\put(4.9,-.05){\line(0,-1){1.5}}
\put(4.8,-.05){\put(-.18,-.18){\tiny$\bullet$}}
\put(5.25,.2){$\tau$}

\thicklines\put(0,-1.5){ 
\qbezier[200](0,0)(-2.2,-.2)(-3.7,.2) \put(-3.8,3){$p(V_r)$}
\qbezier[200](-3.7,.2)(-4.4,.4)(-4.65,1) 
\qbezier[200](-4.65,1)(-4.9,2.2)(-3.7,2.5)
\qbezier[200](-3.7,2.5)(-3.2,2.6)(-2.6,2.2)
\qbezier[200](0,0)(-1.4,1.3)(-2.6,2.2) }
\put(0,-1.5){ 
\qbezier[200](0,0)(-.8,1.5)(-1.1,3.1) 
\qbezier[200](-1.1,3.1)(-1.2,3.7)(-.8,4.2) 
\qbezier[200](-.8,4.2)(0,4.9)(.8,4.2)
\qbezier[200](1.1,3.1)(1.2,3.7)(.8,4.2)
\qbezier[200](0,0)(.8,1.5)(1.1,3.1) }
\put(0,-1.5){ 
\qbezier[200](0,0)(2.2,-.2)(3.7,.2)
\qbezier[50](3.7,.2)(4.25,.45)(4.55,.8)
\put(2.5,2.8){$p(\alpha_r)$}
\qbezier[100](4.8,1.45)(4.5,2.3)(3.7,2.5)
\qbezier[20](4.8,1.45)(4.68,1.125)(4.55,.8)

\qbezier[200](3.7,2.5)(3.2,2.6)(2.6,2.2)
\qbezier[200](0,0)(1.4,1.3)(2.6,2.2) 
} \put(-2.5,-6){\small(a) $T\subset p(\M)$} }

\put(2,0){
\qbezier[300](-2.2,3.2)(0,5.4)(2.2,3.2) 
\qbezier[300](2.2,3.2)(2.8,2)(2.4,0) 
\qbezier[300](2.4,0)(2,-1.4)(1.4,-2.3) 
\qbezier[300](-2.2,3.2)(-2.8,2)(-2.4,0) 
\qbezier[300](-2.4,0)(-2,-1.4)(-1.4,-2.3) 
\qbezier[300](-1.4,-2.3)(0,-4.6)(1.4,-2.3) 
\put(-3,-2){$U_j$}

\put(0,1.8){ \qbezier[80](-.45,.45)(0,.87)(.45,.45)
\qbezier[80](.45,.45)(.7,.16)(.5,-.45)
\qbezier[80](-.45,.45)(-.7,.16)(-.5,-.45)
\qbezier[80](.5,-.45)(.4,-.9)(.2,-1.4)
\qbezier[80](-.5,-.45)(-.4,-.9)(-.2,-1.4)
\qbezier[80](.2,-1.4)(0,-1.8)(-.2,-1.4)

\put(.6,.16){\line(1,0){1.95}} \thicklines
\put(.55,-.45){\line(3,-1){1.95}} \thinlines
\put(.4,-1){\line(3,-2){1.88}}
\put(1.05,-.57){\put(-.17,-.18){\tiny$\bullet$}}\put(1.05,-.57){\line(3,4){1.5}}\put(2.50,1.63){$p(t_j\!)$}
\put(1.5,-1.65){$\tau_j$}

\put(3.2,.16){\line(0,-1){2.5}} \put(3.5,-1.5){$T_j$}
\put(3.2,.135){\line(-1,0){.2}} \put(3.2,-2.34){\line(-1,0){.2}} }

\put(0,-1.5){\put(-.18,-.18){\tiny$\bullet$}} 
\put(-1.2,-2.3){$p(x_0)$} \thicklines
\put(0,-1.5){ 
\qbezier[200](0,0)(-.8,1.5)(-1.1,3.1) 
\qbezier[200](-1.1,3.1)(-1.2,3.7)(-.8,4.2)
\put(-1.3,4.9){$p(\alpha_j)$}
\qbezier[200](-.8,4.2)(0,4.9)(.8,4.2)
\qbezier[200](1.1,3.1)(1.2,3.7)(.8,4.2)
\qbezier[200](0,0)(.8,1.5)(1.1,3.1) } \put(-2.5,-6){\small(b)
$T_j\subset U_j$} }

\put(13.5,0){
\thinlines \put(-2,1){\line(0,1){2}}
           \put(-2,-1){\line(0,-1){2}}
\thicklines \put(2,1){\line(0,-1){2}}
\put(2.3,.4){$\Gamma_j(\tau_j)$}
\qbezier[300](2,0)(.5,-1.5)(-2,-2)
\put(-4.5,-2.4){$\gamma(\alpha_j)$}
\qbezier[3](2,0)(1.7,.3)(1.4,.5) \qbezier[300](1.4,.5)(0,1.6)(-2,2)
\thinlines \qbezier[300](2,1)(.5,-.5)(-2,-1)
\qbezier[300](2,-1)(.5,-2.5)(-2,-3) \qbezier[300](2,1)(0,2.5)(-2,3)
\qbezier[10](2,-1)(1.3,-.4)(.72,0)
\qbezier[200](.72,0)(-.4,.7)(-2,1) \put(-3,-6){\small(c)
$\Gamma_j(T_j)\subset\N$} }
\end{picture}
\end{center}
\caption{The strips $T$, $T_j$ and ``folding'' of
$T_j$ via $\Gamma_j$}\label{figure3}
\end{figure}
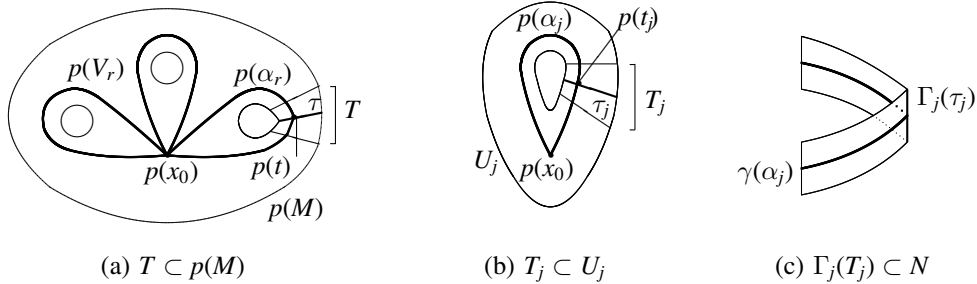

Let $p_\Pi\co  \RR^3\to\Pi$ denote the orthogonal projection.
Without loss of generality, we may assume that the map $p=p_\Pi\circ
i\co  \M\to\Pi$ has the properties~(p1), (p2), (p3) from above.
Moreover,~(p4) holds if $g$ is odd, while the following property
holds if $g$ is even:

(p$4'$)\qua The map $p$ has multiplicity $4$. Moreover, the
set of all points of $p(\M)$, whose preimage under $p$ contains
more than 2 points, lies in a
regular neighbourhood $T$ in $p(\M)$ of a simple arc $\tau\subset
p(\M)$, where the endpoints of $\tau$ lie on the boundary of
$p(\M)$, $\tau$ intersects the graph $p(V_r)$ at the unique point
$p(t)$, for some $t\in\alpha_r\setminus\{x_0\}$, and the
intersection of $\tau$ and $p(\alpha_r)$ at the point $p(t)$ is
transverse; see \fullref{figure3}~(a).

\begin{Pro} \label{pro:Zieschang}
Suppose that $\M$ is an (orientable or nonorientable) closed
surface of genus $g$, and $f\co  \M\to \N$ has absolute degree
$A(f)=0$. Then there exists a self-homeomorphism $\varphi$ of $\M$
and a map $\gamma\co  V_r\to \N$ such that $f\simeq \gamma\circ
\rho\circ\varphi$. Here $r=2g$ if $\M$ is orientable, $r=[\unfrac g2]$
if $\M$ is nonorientable, and $\rho\co  \M\to V_r$ is the
retraction defined above.
\end{Pro}

\Proof Since $A(f)=0$, it follows from~\cite{K} or~\cite{Ep} that
$f$ is homotopic to a map $h$ which is not surjective; thus
$h(\M)\subset \N^*=\N\setminus \IntD^2$ for an appropriate disk
$D^2\subset \N$. Since the fundamental group of $\N^*$ is a free
group, we obtain a homomorphism $h_\#\co   \pi_1(\M) \to
\pi_1(\N^*)$ to the free group $\pi_1(\N^*)$.

Suppose that $\M$ is orientable. It has been proved in Satz~2
of Zieschang \cite{Z} using the Nielsen method (see also Zieschang, Vogt and Coldewey~\cite{ZVC}, or
Proposition~1.2 of~Grigorchuk, Kurchanov and Zieschang \cite{GriKurZie}) that there is a sequence of
``elementary moves'' of the system of generators
$a_1,b_1,\dots,a_g,b_g$ and the corresponding sequence of
``elementary moves'' of the system of cuts
$\alpha_1,\beta_1,\dots,\alpha_g,\beta_g$ on $\M$ (see above), such
that the resulting system of cuts
$\tilde\alpha_1,\tilde\beta_1,\dots,\tilde\alpha_g,\tilde\beta_g$ is
also canonical (this means there exists a self-homeomorphism
$\varphi$ of $\M$ such that $\alpha_j=\varphi(\tilde\alpha_j)$,
$\beta_j=\varphi(\tilde\beta_j)$), and the loops
$\smash{h|_{\tilde\beta_j}}\co \tilde\beta_j \to \N^*$ are contractible in $\N^*$. From this, using
the fact that $\pi_2(\N^*)=0$, one can prove that $h\simeq
\gamma\circ \rho\circ\varphi$ where $\gamma\,:=\,h|_{V_g}$.

Suppose that $\M$ is nonorientable. The method to prove Satz~2
of~\cite{Z} can be successfully applied to construct a canonical
system of cuts $\tilde\alpha_1,\tilde\beta_1,\dots,
\tilde\alpha_{[\unfrac g2]},\tilde\beta_{[\unfrac g2]}$, $\tilde\beta_0$
on $\M$ (this means there exists a homeomorphism $\varphi$ of $\M$
with $\alpha_j=\varphi(\tilde\alpha_j)$,
$\beta_j=\varphi(\tilde\beta_j)$) such that the loops
$\smash{h|_{\tilde\beta_j}}\co \tilde\beta_j \to \N^*$ are contractible in $\N^*$;
see Ol'shanski{\u\i}~\cite{Ol'shanski} or Proposition~1.5 of~\cite{GriKurZie}.
(Again the curve $\beta_0$ is considered only if $g$ is odd.)
Similarly to the orientable case, this implies that $h\simeq
\gamma\circ \rho\circ\varphi$ where $\gamma\,:=\,h|_{V_g}$. \end{proof}

\begin{Thm} \label{thm:A=0}
Suppose that $f\co  \M\to \N$ has absolute degree $A(f)=0$. Then
$2\le \MMR[f]\le 4$.
\end{Thm}

\Proof Suppose that $h$ is homotopic to $f$ and has multiplicity
$1$. Then $h$ is a homeomorphism onto $h(\M)$. It follows from
Brouwer's Theorem on Invariance of Domain~\cite{Brouwer} that $h$ is
surjective and, therefore, it is a homeomorphism. Then $A(h)=1$, a
contradiction. Therefore $\MMR[f]\ge2$.

Let us prove the second inequality. Since $A(f)=0$, by
\fullref{pro:Zieschang}, $f\simeq\gamma\circ
\rho\circ\varphi$ for a self-homeomorphism $\varphi$ of $\M$, the
retraction $\rho\co  \M\to V_r$, and a map $\gamma\co  V_r\to
\N$, where $r=g$ if $\M$ is an orientable surface of genus $g$,
$r=[\unfrac g2]$ if $\M$ is a nonorientable surface of genus $g$.
Without loss of generality, we may assume that $\gamma$ has the
following properties:

$(\gamma1)$\qua There exists a homeomorphism $\psi$ of the
neighbourhood $U$ of $x_0$ in $\M$ onto a neighbourhood of
$\gamma(x_0)$ in $\N$ such that $\gamma|_{V_r\cap U}=\psi|_{V_r\cap
U}$. In other words, $\gamma|_{V_r\cap U}$ extends to an embedding
$\psi\co  U\to \N$;

$(\gamma2)$\qua The restriction of $\gamma$ onto each curve
$\alpha_1,\dots,\alpha_r$ is an immersion $S^1\to\N$. Moreover,
$\gamma$ has multiplicity $\le2$, and it has only finitely many
double points (ie\ pairs of distinct points of $V_r$ having the
same image).

\textbf{Case 1}\qua Suppose that the surface $\M$ is either orientable (thus
$r=g$), or nonorientable with $g$ odd (thus $r=\upnfrac{g-1}2$). In
both cases, the map $p=p_\Pi\circ i\co  \M\to\Pi=\RR^2$ of $\M$
to the plane $\Pi$ has the properties (p1), (p2), (p3), (p4); see
above.

\textbf{Subcase 1}\qua Suppose that $\N$ is orientable. Since every closed curve
$\gamma|_{\alpha_j}$ is orientation-preserving, it follows from the
properties $(\gamma1)$, $(\gamma2)$, (p1), (p3) that the map
$\hat\gamma=\gamma\circ p^{-1}\co  p(V_r)\to \N$ can be extended
to an immersion $\Gamma\co  p(\M)\to \N$ of the regular
neighbourhood $p(\M)$ of $p(V_r)$ in the plane $\Pi$ to $\N$, such
that $\Gamma$ has multiplicity $\le2$.

Consider the composition $\hat\rho=p\circ\rho\co  \M\to\Pi$.
Observe that the maps $\hat\rho$ and $p$ are homotopic as maps
$\M\to p(\M)\subset\Pi$ with the target $p(\M)$, due to
$\hat\rho|_{V_r}=p|_{V_r}$,~(p2), and $\pi_2(p(\M))=0$. From this
and $\gamma=\Gamma\circ p|_{V_r}$, we have
 \begin{equation}\refstepcounter {Thm} \label{eq}
f \simeq \gamma\circ \rho\circ\varphi=\Gamma\circ p\circ
\rho\circ\varphi \simeq \Gamma\circ p\circ\varphi. \tag{\hbox{\bf\theThm}}
 \end{equation}
Since $\varphi$ is bijective and each of $\Gamma$ and $p$ has
multiplicity $\le2$ (see~(p4)), the multiplicity of the composition
$\Gamma\circ p\circ\varphi$ is $\le2\cdot2\cdot1=4$.

\textbf{Subcase 2}\qua Suppose that $\N$ is nonorientable. So in general,
the immersion $\hat\gamma\co  p(V_r)\to \N$ can not be extended
to an immersion of the regular neighbourhood $p(\M)$ of $p(V_r)$ in
$\Pi=\RR^2$. However, due to $(\gamma1)$, $(\gamma2)$, and (p1), we
can extend $\hat\gamma$ to an immersion $\tilde\Gamma\co  p(D\cup
V_r)\to \N$, where $D\subset U$ is a small disk centred at $x_0$.

Now, for each curve $\alpha_j$, we will extend the immersion
$\tilde\Gamma_j=\tilde\Gamma|_{p(D\cup\alpha_j)}\co  p(D\cup\alpha_j)\to\N$
to a regular neighbourhood $U_j\supset p(D)$ of $p(\alpha_j)$ in
$\Pi$ as follows. If the curve $\smash{\gamma|_{\alpha_j}}$ is
orientation-preserving then, similarly to Case~1, the immersion
$\tilde\Gamma_j\co  p(D\cup\alpha_j)\to\N$ can be extended to an
immersion $\Gamma_j\co  U_j\to\N$. If the curve
$\smash{\gamma|_{\alpha_j}}$ is orientation-reversing, let us choose a point
$t_j\in \alpha_j\setminus D$ such that $t_j$ is the only preimage
of the point $\gamma(t_j)$ under $\gamma$. Consider a simple arc
$\tau_j\subset U_j\setminus p(D)$, which transversally intersects
$p(\alpha_j)$ at the only point $p(t_j)$, and whose endpoints lie on
the boundary of $U_j$. Let $T_j$ be a regular neighbourhood of the
arc $\tau_j$ in $U_j\setminus p(D)$, thus $T_j$ is a ``strip'' in
the annulus $U_j$; see \fullref{figure3}~(b). Outside the interior of the strip
$T_j$, we extend $\tilde\Gamma_j$ to an immersion
$\bGamma_j\co   (U_j\setminus T_j)\cup
p(\alpha_j)\to \N$ similarly to above. Now we extend the obtained
immersion $\bGamma_j$ to the whole annulus $U_j$,
giving a map $\Gamma_j\co  U_j\to\N$ which coincides with
$\bGamma_j$ outside $T_j\setminus p(\alpha_j)$ and has
a ``folding'' along the arc $\tau_j\subset T_j$, as shown in
\fullref{figure3}~(c).

Without loss of generality, we may assume that $U_j\subset p(M)$,
and any two annuli $U_j,U_k$ have only the disk $p(D)$ in common.
Since the constructed mappings $\Gamma_j\co  U_j\to\N$ agree on
the common part $p(D)$, they determine an extension
$\bGamma\co  U\to\N$ of the map $\tilde\Gamma$,
where $U=U_1\cup\ldots\cup U_r$ is a regular neighbourhood of
$p(V_r)$ in $\Pi=\RR^2$. The above construction can be performed in
such a way that the map $\bGamma$ has multiplicity
$\le2$, due to~$(\gamma2)$ and the choice of the points
$t_j\in\alpha_j$. Obviously, the map $\bGamma$ can be
extended to the regular neighbourhood $p(\M)$ of $p(D\cup V_r)$
(see~(p3)) and the extended map $\Gamma\co  p(\M)\to\N$ also has
multiplicity $\le2$.

Similarly to Subcase~1, the composition $\Gamma\circ p\circ\varphi$
has multiplicity $\le2\cdot2\cdot1=4$, and \eqref{eq} holds. This
completes the proof in Case~1.

\textbf{Case 2}\qua Suppose that $\M$ is a nonorientable closed surface of even
genus $g$, thus $r=\unfrac g2$, and the map $p=p_\Pi\circ
i\co  \M\to\Pi=\RR^2$ of $\M$ to the plane $\Pi$ has the
properties (p1), (p2), (p3), (p$4'$); see above. We may assume,
without loss of generality, that the map $\gamma\co  V_r\to\N$
has the following additional property:

$(\gamma3)$\qua The point $t\in\alpha_r$ considered in (p$4'$)
is the only preimage of $\gamma(t)$ under $\gamma$, and the
analogous property holds for any point $\tilde t\in \alpha_r\cap
p^{-1}(T)$.

\textbf{Subcase 1}\qua Suppose that $\N$ is orientable. Similarly to Subcase~1
of Case~1, one shows using~$(\gamma1)$, $(\gamma2)$, (p1), (p3) that
the immersion $\hat\gamma=\gamma\circ p^{-1}\co  p(V_r)\to \N$
extends to an immersion $\Gamma\co  p(\M)\to \N$ of multiplicity
$2$, and using~(p2) that~\eqref{eq} holds. Taking into
account~(p$4'$) and $(\gamma3)$, one can show that the multiplicity
of $\Gamma\circ p\circ\varphi$ is $\le4$.

\textbf{Subcase 2}\qua Suppose that $N$ is nonorientable. We proceed as in
Subcase~2 of Case~1. Namely, for those curves $\alpha_j$ whose image
under $\gamma$ is orientation-preserving, we extend the immersion
$\tilde\Gamma_j\co  p(D\cup\alpha_j)\to \N$ to $U_j$, as in
Case~1. For each of the remaining curves $\alpha_j$, we choose a
point $t_j\in\alpha_j\setminus D$ which is the only preimage of
$\gamma(t_j)$ under $\gamma$, and we extend the corresponding
immersion $\tilde\Gamma_j$ to a map
$\bGamma_j\co  U_i\to \N$ having a ``folding'' along
an arc $\tau_j\subset T_j\subset U_j$, which transversally
intersects $p(V_r)$ at the unique point $p(t_j)$; see Case~1. As
above, this allows one to construct a map $\Gamma\co  p(\M)\to\N$
of multiplicity $\le2$ which is an extension of $\hat\gamma$, and to
show that~\eqref{eq} holds. Observe now that, if the curve
$\gamma|_{\alpha_r}$ is orientation-reversing, we can choose the
point $t_r\in\alpha_r$ in such a way that it is ``far enough'' from
the point $t\in\alpha_r$ considered in~(p$4'$). This, together with
$(\gamma3)$, shows that the above construction can be performed in
such a way that the composition $\Gamma\circ p\circ\varphi$ has
multiplicity $\le4$. This completes the proof of 
\fullref{thm:A=0}. \end{proof}

\bibliographystyle{gtart}
\bibliography{link}

\end{document}